\documentclass[reqno]{amsart}
\usepackage{float}
\usepackage{amsmath}
\usepackage{graphicx}
\usepackage{latexsym}
\usepackage{amsfonts}
\usepackage{amssymb}
\theoremstyle{plain}
\newtheorem{theorem}{Theorem}
\newtheorem{corollary}[theorem]{Corollary}

\theoremstyle{definition}

\newtheorem{definition}[theorem]{Definition}

\newtheorem*{remark*}{Remark}

\begin{document}
\title[Martingale approach to subexponential asymptotics]{Martingale approach to subexponential asymptotics for random walks}
\author[Denisov]{Denis Denisov}
\address{School of Mathematics, Cardiff University, Senghennydd Road
CARDIFF, Wales, UK.
CF24 4AG Cardiff }
\email{DenisovD@cf.ac.uk}

\author[Wachtel]{Vitali Wachtel}
\address{Mathematical Institute, University of Munich, Theresienstrasse 39, D--80333
Munich, Germany}
\email{wachtel@mathematik.uni-muenchen.de}

\begin{abstract}
Consider the random walk $S_n=\xi_1+\cdots+\xi_n$ with 
independent and identically distributed increments and negative mean 
$\mathbf E\xi=-m<0$. Let $M=\sup_{0\le i} S_i$ be the supremum of the random walk.  In
this note we present derivation of asymptotics for $\mathbf P(M>x), x\to\infty$ for  
long-tailed distributions. This derivation is based on the martingale arguments and 
does not require any prior knowledge of the theory of long-tailed distributions. 
In addition the same approach allows to obtain asymptotics for $\mathbf P(M_\tau>x)$, where 
$M_\tau=\max_{0\le i<\tau}S_i$ and $\tau=\min\{n\ge 1: S_n\le 0 \}$. 
\end{abstract}


\maketitle
{\scriptsize
}

\section{Introduction, statement of results and discussion}

Let $\xi,\xi_1$, $\xi_2$, \ldots\ be independent random variables with
a common distribution function~$F$ and negative mean, i.e., 
$\mathbf E\xi=-a<0$. Let $S_n$ denote the random walk with the increments 
$\xi_k$, that is,
$$
S_0=0,\quad S_n=\xi_1+\xi_2+\cdots+\xi_n,\ \ n\geq1.
$$
It follows from the assumption $\mathbf{E}\xi<0$ that the total maximum
$M:=\sup_{n\geq0} S_n$ is finite almost surely. The asymptotic behaviour
of $\mathbf{P}(M>x)$ has been considered by many authors. The first
results are due to Cramer and Lundberg: if there exists $h_0>0$ such
that $\mathbf{E} e^{h_0\xi}=1$ and $\mathbf{E}\xi e^{h_0\xi}<\infty$ then
\begin{equation}
\label{Cramer-Lundberg}
\mathbf{P}(M>x)\sim c_0e^{-h_0x}\quad\text{as }x\to\infty
\end{equation}
for some $c_0\in(0,1)$ and, furthermore,
\begin{equation}
\label{CL-inequality}
\mathbf{P}(M>x)\leq e^{-h_0x}\quad\text{for all }x>0.
\end{equation}
The proof of these statements is based on the following observation: The
assumption $\mathbf{E}e^{h_0\xi}=1$ implies that the sequence $e^{h_0S_n}$
is a martingale. Applying the Doob inequality we obtain immediately
(\ref{CL-inequality}). The same martingale property allows one to make an
exponential change of measure, which is used in the proof of 
(\ref{Cramer-Lundberg}).

If the distribution of $\xi$ is long-tailed, i.e.,
$\mathbf{E}e^{h\xi}=\infty$ for all $h>0$, then one can investigate 
$\mathbf{P}(M>x)$ under some additional regularity restrictions on the
tail function $\overline{F}(x):=1-F(x)$.
One of the most popular regularity assumption is the so-called 
subexponentiality of the distribution tails.
\begin{definition}
The distribution function $F$ on ${\bf R}_+$ is called subexponential if
$$
\int_0^x \overline F(x-y)dF(y)\sim 2\overline F(x) 
\quad\mbox{ as }x\to\infty.
$$
\end{definition}
The following result is known in the literature as Veraverbecke's theorem:
Let $F_I$ be defined by the tail 
$\overline{F_I}(x):=\min\left(1,\int_x^\infty\overline{F}(y)dy\right)$, $x>0$. 
If $F_I$ is subexponential then
\begin{equation}
\label{Veraverbecke}
\mathbf P(M>x)\sim\frac{1}{a}\overline{F_I}(x)\quad\text{as }x\to\infty.
\end{equation}

We next turn to the maximum of the positive excursion of the random walk.
Let
$$
\tau:=\inf\{n\ge 1: S_n\le 0\}
$$
and 
$$
M_\tau:=\max_{0\le n<\tau}S_n. 
$$
If the Cramer-Lundberg condition holds then one can derive the asymptotics
for $\mathbf{P}(M_\tau>x)$ from that for the total maximum $M$. This way
has been suggested first by Iglehart \cite{Igl}. Namely, it follows from the Markov
property that
$$
\mathbf{P}(M>x)=\mathbf{P}(M_\tau>x)+
\int_{-\infty}^0\mathbf{P}(M>x-y)\mathbf{P}(S_\tau\in dy,M_\tau\leq x).
$$
Thus,
$$
\mathbf{P}(M_\tau>x)=\mathbf{P}(M>x)\left(1-
\int_{-\infty}^0\frac{\mathbf{P}(M>x-y)}{\mathbf{P}(M>x)}\mathbf{P}(S_\tau\in dy,M_\tau\leq x)\right).
$$
Noting that (\ref{Cramer-Lundberg}) yields
$$
\lim_{x\to\infty}\frac{\mathbf{P}(M>x-y)}{\mathbf{P}(M>x)}=e^{h_0y}\quad
\text{for every }y<0,
$$
and applying the dominated convergence, we obtain
$$
\int_{-\infty}^0\frac{\mathbf{P}(M>x-y)}{\mathbf{P}(M>x)}\mathbf{P}(S_\tau\in dy,M_\tau\leq x)
\sim\mathbf{E} e^{h_0S_\tau}.
$$
As a result we get
\begin{equation}
\label{CL-local}
\mathbf{P}(M_\tau>x)\sim\left(1-\mathbf{E} e^{h_0S_\tau}\right)\mathbf{P}(M>x)
\sim \left(1-\mathbf{E} e^{h_0S_\tau}\right)c_0 e^{-h_0x}.
\end{equation}

It turns out that Iglehart's approach can not be applied to heavy-tailed random
walks without further restrictions on the distribution of $\xi$. Here one has to
assume that $F$ is {\it strong subexponential}. This class of distribution 
functions was introduced by Kl\"uppelberg \cite{Kl}.

\begin{definition}
The distribution function $F$ on ${\bf R}$ belongs to the
class ${\mathcal S}\,^*$  if
\begin{eqnarray}
\label{s_sub}
\int_0^x \overline F(x-y)\overline F(y)dy
&\sim& 2a_+ \overline F(x) 
\quad\mbox{ as }x\to\infty,
\end{eqnarray}
where $a_+=\int_0^\infty \overline F(y)dy\in(0,\infty).$ 
\end{definition}

Denisov \cite{D} adopted Iglehart's reduction from $M_\tau$ to $M$ to
the class of strong subexponential distributions: If $F\in{\mathcal S}\,^*$
then
\begin{equation}\label{eq_asm}
\mathbf P(M_\tau>x)\sim\mathbf E\tau \overline F(x),\quad x\to\infty. 
\end{equation}
The asymptotics (\ref{eq_asm}) were found first by Asmussen \cite{A1} 
for $F\in\mathcal S^*$ and by Heath, Resnick and Samorodnitsky \cite{HRS} for 
regularly varying $F$. An extension of this result to the general stopping 
time can be found in Foss and Zachary \cite{FZ}, and in Foss, Palmowsky and 
Zachary \cite{FPZ}. These extensions rely on (\ref{eq_asm}). 

The {\it main purpose} of the present note is to give alternative proofs 
of (\ref{Veraverbecke}) and (\ref{eq_asm}) using martingale techniques.

In order to state our main result we introduce some notation.
For any $y>0$ let  
$$
\mu_y:=\min\{n\ge 0: S_n>y\}. 
$$
The latter stopping time is naturally connected with the supremum since 
$$
\mathbf P(M>x)=\mathbf P(\mu_x<\infty). 
$$

Let 
$$
\overline F_s(x):=\int_x^\infty\overline F(u)du
$$
and 
\begin{equation}
\label{defG}
G_{c}(x)= \begin{cases} 
                       \overline F_s(x), & \mbox{if } x\ge 0 \\
                        c, & \mbox{if } x<0 
          \end{cases}. 
\end{equation}
Define also 
\begin{equation}
\label{defG1}
\widehat{G}_c(x):=\min\{G_c(x),c\}.
\end{equation}
\begin{theorem}
\label{main}
Assume that $F$ is long-tailed. For any $\varepsilon>0$ 
there exists $R>0$ such that the stopped sequence 
\begin{equation}\label{sub}
\widehat{G}_{a+\varepsilon}(x-S_{n\wedge \mu_{x-R}}) \quad \text{is a submartingale}. 
\end{equation}
Assume in addition  that $F\in\mathcal S^*$.  
For any $\varepsilon>0$ 
there exists $R>0$ such that the stopped sequence 
\begin{equation}
\label{super}
\widehat G_{a-\varepsilon}(x-S_{n\wedge \mu_{x-R}}) \quad \text{is a supermartingale}. 
\end{equation}
\end{theorem}
Having constructed super- and submartingale we can obtain subexponential 
asymptotics for $M$ and $M_\tau$ by applying the optional stopping theorem.
\begin{corollary}
\label{cor.bounds}
For any long-tailed distribution function $F$ with negative mean,
\begin{equation}
\label{lower}
\liminf_{x\to\infty} \frac{\mathbf P(M>x)}{\overline F_s(x)}\ge \frac{1}{a}. 
\end{equation}
Assume in addition that  $F\in\mathcal S^*$. 
Then,
$$
\mathbf P(M>x)\sim\frac{1}{a}\overline F_s(x),\quad x\to\infty.
$$
\end{corollary}
To the best of our knowledge, all existing in the literature proofs
of the Veraverbecke theorem are based on representations via geometric sums.
More precisely, $\mathbf{P}(M>x)$ can be estimated by
$\sum_{n=1}^\infty (1-p)p^{n}\mathbf{P}(Y_1+Y_2+\ldots+Y_n>x)$, where
$p\in(0,1)$ and $Y_i$ are independent identically distributed random
variables with $\mathbf{P}(Y_1>x)\sim \overline{F_s}(x)$. In order to
obtain (\ref{Veraverbecke}) from that geometric sum one uses the following two
properties of subexponential distributions:
\begin{itemize}
\item[(a)]$\mathbf{P}((Y_1+Y_2+\ldots+Y_n>x))\sim n\mathbf{P}(Y_1>x)$ for every fixed $k$,
\item[(b)]For every $\varepsilon>0$ there exists $C(\varepsilon)<\infty$ such that
$$
\mathbf{P}((Y_1+Y_2+\ldots+Y_n>x))\leq C(\varepsilon)(1+\varepsilon)^n\mathbf{P}(Y_1>x).
$$
\end{itemize}
A recent elegant proof based on (a) and (b) can be found in \cite{Z}.
Our proof does not use any property of $F$ besides (\ref{s_sub}).

Unfortunately, our method does not allow us to derive (\ref{Veraverbecke}) for
the whole class of subexponential distributions. The condition $F_I\in \mathcal S$ 
and $F_I\in \mathcal S^\ast$ are close but do not coincide, see Section~6 in \cite{DFK}.
But we can apply the same construction
to $M_\tau$ and, as it is known in the literature, the strong subexponentiality
is optimal for asymptotics (\ref{eq_asm}).
\begin{corollary} \label{res}
Let $F\in \mathcal {\mathcal S}^*$. Then 
\begin{equation}\label{asm}
\mathbf P(M_\tau>x)\sim\mathbf E\tau \overline F(x).
\end{equation}
\end{corollary}
It is worth mentioning that, in contrast to all previous proofs,
our approach to (\ref{asm}) is direct, i.e., it
does not use any knowledge on the asymptotic behaviour of $M$.

One of the important advantages of the martingale approach is the possibility
to obtain non-asymptotic inequalities for $\mathbf{P}(M>x)$ and $\mathbf{P}(M_\tau>x)$.
For example, it follows from (\ref{sub}) that for every $\varepsilon>0$ there exists
$R>0$ such that (see the proof of Corollary \ref{cor.bounds})
\begin{equation}
\label{lower1}
\mathbf{P}(M>x)\geq\frac{\overline{F_s}(x+R)}{a+\varepsilon},\quad x>0.
\end{equation}
Using a supermartingale property of $G_{a-\varepsilon}$ we obtain the following
upper bound
\begin{equation}
\label{upper1}
\mathbf{P}(M>x)\leq\frac{\overline{F_s}(x-R')}{a-\varepsilon},\quad x>R'.
\end{equation}
Of course, in order to apply these inequalities, one has to know how to compute
$R$ and $R'$ for given values of $\varepsilon$. And we believe that one can do it 
rather easy for certain subclasses of $\mathcal{S}^*$, e.g., for regularly varying 
or Weibull tails.

Foss, Korshunov and Zachary \cite{FKZ} have shown that the inequality
$$
\mathbf{P}(M>x)\geq\frac{\overline{F_s}(x)}{a+\overline{F_s}(x)},\quad x>0
$$
holds without any restriction on the distribution function $F$, see Theorem 5.1 in \cite{FKZ}. 
This bound is better
than (\ref{lower1}). It's proof is based on the fact, that the distribution of $M$ is the
stationary distribution of the Lindley recursion $W_{n+1}=(W_n+\xi_{n+1})^+$. This property
of $M$ can be written as follows: Let $\xi'$ a copy of $\xi$, which is independent of $M$.
Then $\mathcal{L}(M)=\mathcal{L}((M+\xi'))$. This can be seen as a martingale property:
Define $\pi(x):=\mathbf{P}(M>x)$. Then the sequence $\pi(x-S_{n\wedge\mu_x})$ is a martingale.

Using (\ref{super}), one gets for all $x>R'$ the inequality
$$
\mathbf{P}(M_\tau>x)\leq\frac{\overline{F_s}(x-R')-\mathbf{E}\overline{F_s}(x-R'-S_\tau)}{a-\varepsilon}.
$$
And an upper estimate for the difference in the nominator is easy to get:
\begin{align*}
\overline{F_s}(x-R')-\mathbf{E}\overline{F_s}(x-R'-S_\tau)=
\mathbf{E}\left[\int_{x-R'}^{x-R'-S_\tau}\overline{F}(z)dz\right]
\leq \overline{F}(x-R')\mathbf{E}[-S_\tau].
\end{align*}
Applying the Wald identity, we obtain
\begin{equation}
\mathbf{P}(M_\tau>x)\leq\frac{a}{a-\varepsilon}\mathbf{E}\tau\overline{F}(x-R').
\label{upper2}
\end{equation}
A lower bound is not as obvious. Here we can conclude from (\ref{sub}) that
\begin{equation}
\label{lower2}
\mathbf{P}(M_\tau>x)\geq\frac{\overline{F_s}(x+R)-\mathbf{E}\overline{F_s}(x+R-S_\tau)}{a+\varepsilon}.
\end{equation}
Thus one needs an appropriate estimate for the difference in the nominator. 

Martingale approach has been used also by Kugler and Wachtel \cite{KW} in deriving upper bounds
for $\mathbf{P}(M>x)$ and $\mathbf{P}(M_{\tau_z}>x)$, where $\tau_z:=\min\{k:S_n\leq -z\}$ under
the assumption that some power moments of $\xi$ are finite. Their strategy is completely different: 
They truncate the summands $\xi_i$ in order to construct an exponential supermartingale for the 
random walk with truncated increments.
\section{Proofs.}
\subsection{Proof of Theorem \ref{main}.}
Fix $\varepsilon>0$. To prove the submartingale property we need to show that 
\begin{equation}
  \label{sub1}
\mathbf E \widehat{G}_{a+\varepsilon}(x-y-\xi)\ge \widehat{G}_{a+\varepsilon}(x-y)
\end{equation}
for all $y\leq x-R$.

Put, for brevity, $t:=x-y\ge R$. By  the definition (\ref{defG}),  
\begin{align*}
\mathbf E \widehat{G}_{a+\varepsilon}(t-\xi)&=(a+\varepsilon)\mathbf P(\xi>t-r_c)
+\int_{-\infty}^{t-r_c}F(dz)\overline{F_s}(t-z)\\
&=(a+\varepsilon)\overline F(t-r_c)+\left(\int_0^{t-r_c}+\int_{-\infty}^0\right)
F(dz)\overline{F_s}(t-z),
\end{align*}
where $r_c:=\min\{x\geq0:\overline{F_s}(x)\leq c\}$.
Integrating the first integral by parts, we obtain 
\begin{align*}
\int_{0}^{t-r_c}F(dz)\overline{F_s}(t-z)=
\overline F(0)\overline {F_s}(t-r_c)-\overline F(t-r_c)\overline{F_s}(0)
+\int_0^{t-r_c}dz \overline F(z)\overline F(t-z).
\end{align*}
Integrating the second integral by parts, we obtain 
\begin{align*}
  \int_{-\infty}^{0}F(dz)\overline{F_s}(t-z)=F(0)\overline{F_s}(t)-\int_{-\infty}^0dz\overline F(t-z)F(z).
\end{align*}
Combining the above inequalities, we get
\begin{align}
\mathbf E \widehat{G}_{a+\varepsilon}(t-\xi)
&=(a+\varepsilon)\overline F(t-r_c)-\overline F(t-r_c)\overline {F_s}(0)
+\overline{F}(0)\overline {F_s}(t-r_c)\nonumber\\
\label{after.by.parts}
&+\int_0^{t-r_c}dz \overline F(z)\overline F(t-z)
+F(0)\overline{F_s}(t)-\int_{-\infty}^0dz\overline F(t-z)F(z).
\end{align}
It is clear that
$$
\int_{-\infty}^0dz\overline F(t-z)F(z)\leq \overline F(t)\int_{-\infty}^0dzF(z)=a_-\overline{F}(t).
$$
Further,
$$
\overline{F}(0)\overline {F_s}(t-r_c)+F(0)\overline{F_s}(t)=
\overline{F_s}(t)+\overline{F}(0)\int_{t-r_c}^t\overline{F}(z)dz
$$
and
$$
\int_0^{t-r_c}dz \overline F(z)\overline F(t-z)=
\int_0^{t}dz \overline F(z)\overline F(t-z)-
\int_{t-r_c}^tdz \overline F(z)\overline F(t-z).
$$
Now, put $a_+:=\overline{F_s}(0), a_-:=\int_{-\infty}^0dzF(z)$ and note that $a=a_--a_+.$ 
Consequently,
\begin{align*}
\mathbf E \widehat{G}_{a+\varepsilon}(t-\xi)
&\ge \overline{F_s}(t)+(a+\varepsilon)\overline F(t-r_c)-a_+\overline F(t-r_c)-a_-\overline F(t)\\
&\hspace{1cm}+2\int_0^{t/2}dz \overline F(z)\overline F(t-z)+
\int_{t-r_c}^t\overline{F}(z)\left(\overline{F}(0)-\overline F(t-z)\right)dz\\
&\ge \overline{F_s}(t)+(-2a_++\varepsilon)\overline F(t-r_c)+
2\overline F(t)\int_0^{t/2}dz\overline{F}(z).
\end{align*}
Now, taking $R_1$ sufficiently large, we can ensure that 
$$
2\int_0^{t/2}\overline{F}(z)dz\geq 2a_+-\frac{\varepsilon}{2}\quad\text{for all }t\geq R_1.
$$ 
Furthermore, we can choose $R_2$ so large that
$$
\frac{\overline{F}(t-r_c)-\overline{F}(t)}{\overline{F}(t)}\leq\frac{\varepsilon}{4a_+}.
$$
As a result, for $t>\max\{R_1,R_2\}$ we have
$$
\mathbf E \widehat{G}_{a+\varepsilon}(t-\xi)\geq \overline{F_s}(t)
$$
This proves (\ref{sub}).

To prove the supermartingale property 
it sufficient to show  that 
\begin{equation}
  \label{super1}
\mathbf E G_{a-\varepsilon}(x-y-\xi)\le G_{a-\varepsilon}(x-y)
\end{equation}
for all $y\leq x-R$.
Using (\ref{after.by.parts}) with $r_c=0$, we obtain
\begin{align*}
\mathbf E G_{a-\varepsilon}(t-\xi)
&=G_{a-\varepsilon}(t)+(a-\varepsilon-a_+)\overline F(t)\\
&+\int_0^{t}dz \overline F(z)\overline F(t-z)
-\int_{-\infty}^0dz\overline F(t-z)F(z).
\end{align*}
According to the definition of $\mathcal S^*$ there exists $R_1$ such that
$$
\int_0^{t}dz \overline F(z)\overline F(t-z)\le (2a_++\varepsilon/2)\overline F(t)
$$
for all $t\geq R_1$.
Furthermore, since $F$ is long-tailed, we have
$$
\lim_{t\to\infty}\frac{1}{\overline{F}(t)}\int_{-\infty}^0dz\overline F(t-z)F(z)=
\int_{-\infty}^0dzF(z)=a_-.
$$
Therefore, there exists $R_2$ such that
$$
\int_{-\infty}^0dz\overline F(t-z)F(z)\ge (a_-+\varepsilon/2)\overline F(t),
\quad t\geq R_2.
$$
This immediately implies (\ref{super}) with $R=\max\{R_1,R_2\}$.

\subsection{Proof of Corollary~\ref{cor.bounds}.} 
Fix $\varepsilon>0$ and pick $R$ such that 
$$
Y_n=\widehat{G}_{a+\varepsilon}(x-S_{n\wedge \mu_{x-R}})
$$
is a submartingale. Then,
\begin{align*}
\overline{F_s}(x)=\widehat{G}_{a+\varepsilon}(x) &=  \mathbf EY_0\le \mathbf EY_{\infty}\\
&=\mathbf{E}\left[\widehat{G}_{a+\varepsilon}(x-S_{\mu_{x-R}}),\mu_{x-R}<\infty\right]\\
&\le (a+\varepsilon)\mathbf P(\mu_{x-R}<\infty).
\end{align*}
Hence, 
$$
\mathbf P(M>x)=\mathbf P(\mu_x<\infty)\ge\frac{1}{a+\varepsilon}\overline{F_s}(x+R).
$$
Letting $x$ to infinity we obtain, 
$$
\liminf_{x\to\infty}\frac{\mathbf P(M>x)}{\overline{F_s}(x)}\ge a+\varepsilon.
$$
Since $\varepsilon>0$ is arbitrary the lower bound in (\ref{lower}) holds. 

To prove the corresponding upper bound  
fix $\varepsilon>0$ and pick $R$ such that 
$Y_n=G_{a-\varepsilon}(x-S_{n\wedge \mu_{x-R}})$ is a supermartingale.  Then,
\begin{align}
\nonumber
\overline{F_s}(x)=G_{a-\varepsilon}(x) &=  \mathbf EY_0\ge \mathbf EY_{\infty}\\
\nonumber
&=(a-\varepsilon)\mathbf P(\mu_{x-R}<\infty, S_{\mu_{x-R}}>x)\\
\nonumber
&\hspace{0.5cm}+\mathbf E\left[ \overline{F_s}(x-S_{\mu_{x-R}});\mu_{x-R}<\infty, S_{\mu_{x-R}}\in(x-R,x] \right]\\
\nonumber
&\ge (a-\varepsilon)\mathbf P(\mu_{x-R}<\infty, S_{\mu_{x-R}}>x)\\
\label{upper}
&\hspace{0.5cm}+\overline{F_s}(R)\mathbf P(\mu_{x-R}<\infty, S_{\mu_{x-R}}\in(x-R,x]).
\end{align}
Let $r>0$ be a number which we pick later. Then, 
\begin{align*}
\mathbf P(M>x+r)&\le 
\mathbf P(\mu_{x-R}<\infty, S_{\mu_{x-R}}>x)\\
&\hspace{1cm}+\mathbf P(\mu_{x-R}<\infty, S_{\mu_{x-R}}\in(x-R,x],M>x+r)\\
&\le 
\mathbf P(\mu_{x-R}<\infty, S_{\mu_{x-R}}>x)\\
&\hspace{1cm}+
\mathbf P(M>r)\mathbf P(\mu_{x-R}<\infty, S_{\mu_{x-R}}\in(x-R,x]),
\end{align*}
where we use the strong Markov property. 
Now pick sufficiently large $r$ such that $\mathbf P(M>r)\le \overline{F_s}(R)/(a-\varepsilon)$. Then,
\begin{align*}
 \mathbf P(M>x+r)&\le 
\mathbf P(\mu_{x-R}<\infty, S_{\mu_{x-R}}>x)\\
&\hspace{1cm}+
\frac{\overline{F_s}(R)}{a-\varepsilon}\mathbf P(\mu_{x-R}<\infty, S_{\mu_{x-R}}\in(x-R,x]).\\
\end{align*}
Combining this with (\ref{upper}), we get
$$
\mathbf P(M>x+r)\leq\frac{\overline{F_s}(x)}{a-\varepsilon}.
$$
Letting $x$ to infinity we obtain, 
$$
\limsup_{x\to\infty}\frac{\mathbf P(M>x)}{\overline{F_s}(x)}\le \frac{1}{a-\varepsilon}.
$$
Since $\varepsilon>0$ is arbitrary the upper bound holds. 

\subsection{Proof of Corollary~\ref{res}}
We start with a lower bound. 
Fix $\varepsilon>0$ and pick $R$ such that 
$Y_n=\widehat{G}_{a+\varepsilon}(x-S_{n\wedge \mu_{x-R}})$ is a submartingale.  Then,
\begin{align*}
\overline{F_s}(x)&=\widehat{G}_{a+\varepsilon}(x)=  \mathbf EY_0\le \mathbf EY_{\tau}\\
&\leq(a+\varepsilon)\mathbf P(\mu_{x-R}<\tau)+\mathbf E\overline{F_s}(x-S_\tau).
\end{align*}
Hence,
\begin{align*}
 \mathbf P(M_\tau>x+R)=\mathbf P(\mu_{x+R}<\tau)\ge 
\frac{1}{a+\varepsilon}\left(\overline{F_s}(x)-\mathbf E\overline{F_s}(x-S_\tau)\right).
\end{align*}
Now 
\begin{align}
\overline{F_s}(x)-\mathbf E\overline{F_s}(x-S_\tau)
&=\int_0^{\infty}\mathbf P(S_\tau\in -dt)\left(\overline{F_s}(x)-\overline{F_s}(x+t)\right)\nonumber\\
&\sim |\mathbf  ES_\tau|\overline F(x),\quad x\to\infty.
\label{eq0}
\end{align}
By the Wald's identity $|\mathbf ES_\tau|=a\mathbf E\tau$. Therefore, 
$$
\liminf_{x\to \infty}\frac{\mathbf P(M_\tau>x)}{\overline F(x)}\ge\frac{a\mathbf E\tau}{a+\varepsilon}.
$$
Since $\varepsilon>0$ is arbitrary we obtain the lower bound 
$$
\liminf_{x\to \infty}\frac{\mathbf P(M_\tau>x)}{\overline F(x)}\ge\mathbf E\tau.
$$

To show the upper bound fix $\varepsilon>0$ and pick $R$ such that 
$Y_n=G_{a-\varepsilon}(x-S_{n\wedge \mu_{x-R}})$ is a supermartingale.  
Then,
\begin{align*}
\overline{F_s}(x)&=G_{a-\varepsilon}(0)=  \mathbf EY_0\ge \mathbf EY_{\tau}\\
&=(a-\varepsilon)\mathbf P(\mu_{x-R}<\tau, S_{\mu_{x-R}}>x)\\
&\hspace{0.5cm}+\mathbf E\left[ \overline{F_s}(x-S_{\mu_{x-R}});\mu_{x-R}<\tau, S_{\mu_{x-R}}\in(x-R,x] \right]\\
&\hspace{0.5cm}+\mathbf E\overline{F_s}(x-S_\tau)\\
&\geq(a-\varepsilon)\mathbf P(\mu_{x-R}<\tau, S_{\mu_{x-R}}>x)+\mathbf E\overline{F_s}(x-S_\tau)\\
&\hspace{0.5cm}+\overline{F_s}(R)\mathbf P(\mu_{x-R}<\tau, S_{\mu_{x-R}}\in(x-R,x]).
\end{align*}
Similarly to the corresponding argument in the proof of Corollary~\ref{cor.bounds},
\begin{align*}
\mathbf P(M_\tau>x+r)&\le 
\mathbf P(\mu_{x-R}<\tau, S_{\mu_{x-R}}>x)\\
&\hspace{1cm}+\mathbf P(\mu_{x-R}<\tau, S_{\mu_{x-R}}\in(x-R,x],M_\tau>x+r)\\
&\le 
\mathbf P(\mu_{x-R}<\tau, S_{\mu_{x-R}}>x)\\
&\hspace{1cm}+
\mathbf P(M>r)\mathbf P(\mu_{x-R}<\tau, S_{\mu_{x-R}}\in(x-R,x]),
\end{align*}
Consequently, 
$$
\mathbf P(M_\tau>x+r)\le \frac{1}{a+\varepsilon}\left(\overline{F_s}(x)-\mathbf E\overline{F_s}(x-S_\tau)\right).
$$
Now, we can apply (\ref{eq0}) and obtain
$$
\limsup_{x\to \infty}\frac{\mathbf P(M_\tau>x)}{\overline F(x)}\le\frac{a\mathbf E\tau}{a-\varepsilon}.
$$
Since $\varepsilon>0$ is arbitrary we obtain the upper bound 
$$
\limsup_{x\to \infty}\frac{\mathbf P(M_\tau>x)}{\overline F(x)}\le\mathbf E\tau.
$$

\qed

{\bf Acknowledgements.}
This work was partially supported by DFG. 
This work was done while Denis Denisov was visiting the University of Muenchen.

\end{document}